  \documentclass{gretsi}
 \usepackage[utf8]{inputenc} 
 \usepackage[english,french]{babel}   
  \usepackage{times}			
  
  \usepackage{amsmath,graphicx}

\usepackage{dsfont}
\usepackage{amssymb}
\usepackage{amsbsy}
 \usepackage{algorithmic} 
\usepackage{algorithm,caption}
 \usepackage{url}
\usepackage{scalefnt}

\floatname{algorithm}{Algorithme}

\newcommand{\algorithmicoutput}{\textbf{Sortie:}}
\newcommand{\OUTPUT}{\item[\algorithmicoutput]}

\newcommand{\algorithmicinput}{\textbf{Entrées:}}
\newcommand{\INPUT}{\item[\algorithmicinput]}

  \def\x{{\mathbf x}}
\def\y{{\mathbf y}}
\def\z{{\mathbf z}}
\def\ub{{\mathbf u}}

\def\h{{\mathbf h}}
\def\V{{\mathbf V}}

\def\R{{\mathbb R}}

\def\xib{{\boldsymbol{\xi}}}
 
 \usepackage[T1]{fontenc}

\titre{Modélisations de textures par champ gaussien\\ à orientation locale prescrite}

\auteur{\coord{Kévin}{Polisano}{1},
        \coord{Marianne}{Clausel}{1},
    \coord{Valérie}{Perrier}{1},
    \coord{Laurent}{Condat}{2}}
    
\adresse{\affil{1}{Université Grenoble Alpes et CNRS, Laboratoire Jean Kuntzmann, F-38000, Grenoble}
         \affil{2}{Université Grenoble Alpes et CNRS, GIPSA-lab, F-38000, Grenoble}}

\email{\{Kevin.Polisano,Marianne.Clausel,Valerie.Perrier\}@imag.fr}

\resumefrancais{Nous présentons dans ce papier deux nouveaux modèles de texture orientée, basés sur une nouvelle classe de champs gaussiens, appelés \textit{champs browniens fractionnaires localement anisotropes}, à orientation locale prescrite en chaque point. Ces champs aléatoires sont une version locale d'une classe spécifique de champs anisotropes autosimilaires à incréments stationnaires. La simulation de telles textures s'appuie sur un nouvel algorithme couplant l'utilisation de la notion de champs tangents à la méthode de Cholesky ou plus récemment la méthode des bandes tournantes, cette dernière ayant prouvé son efficacité à générer des textures stationnaires anisotropes. Des applications numériques illustrent la faculté de la méthode à synthétiser des textures à orientation locale prescrite.}

\resumeanglais{This paper presents two new models of oriented texture, based on a new class of Gaussian fields, called {\it locally anisotropic fractional Brownian fields}, with prescribed local orientation at any point. These fields are a local version of a specific class of anisotropic self-similar Gaussian fields with stationary increments. The simulation of such textures is obtained using a new algorithm mixing the tangent field formulation with the Cholesky method or the turning band method, this latter method having proved its efficiency for generating stationary anisotropic textures. Numerical experiments show the ability of the method for synthesis of textures with prescribed local orientation.}

\begin{document}
\maketitle

\section{Textures et champs gaussiens autosimilaires anisotropes}
La modélisation de texture est un problème difficile du traitement d'image. Il existe une variété de méthodes associées à différents champs des mathématiques : des méthodes dites structurelles, statistiques, ou encore à base de patchs, etc. Nous nous intéressons ici aux textures aléatoires présentant un comportement fractal, pour lesquelles les même motifs se retrouvent à différentes échelles, comme c'est le cas pour les nuages ou les reliefs montagneux. L'outil mathématique derrière cette propriété d'autosimilarité est le mouvement brownien fractionnaire, modèle stochastique bien connu introduit dans \cite{mandelbrot1968}. Néanmoins, l'utilisation des champs browniens fractionnaires, qui sont isotropes par définition, n'est pas adaptée à des données présentant de l'anisotropie. D'autres modèles ont été introduits dans la littérature pour y remédier, citons notamment le drap brownien fractionnaire \cite{kamont1996,ayache2002} et le champ brownien fractionnaire anisotrope \cite{Bonami2003}, satisfaisant une propriété d'anisotropie globale, ou dans un autre registre les textures localement parallèles \cite{maurel2011}. La conception mathématique de modèles de textures anisotropes et leur synthèse, sont des éléments essentiels à l'estimation des caractéristiques anisotropes des images. Nous avons introduit dans \cite{polisano:2014} un nouveau modèle. Nous proposons dans cet article deux simulations par champs tangents de ce modèle, dont un est extrait de \cite{polisano:2014}. La brique de base de ce travail s'appuie sur le champ brownien fractionnaire dont nous rappelons la définition et les propriétés.


Le \textit{champ brownien fractionnaire} (CBF) d'indice de Hurst $H\in (0,1)$, désigné par
$B^H$, est l'unique champ gaussien à valeurs réelles vérifiant :\\
 -- presque surement $B^H(0)=0$,\\
 -- incréments stationnaires : pour tout $\z\in \R^2,$\\ $B^H(\cdot + \z)-B^H(\z)\overset{\mathcal{L}}{=} B^H(\cdot)-B^H(0),$\\
 -- autosimilarité d'ordre $H$ : $\forall \lambda \in \R^{\star}, B^H(\lambda \cdot)\overset{\mathcal{L}}{=} \lambda^{H}B^H(\cdot),$\\
 -- isotropie : pour toute rotation $R$ de $\R^2$, $B^H\circ R\overset{\mathcal{L}}{=} B^H,$\\
où $\overset{\mathcal{L}}{=}$ désigne l'égalité en loi. \\
D'après \cite{samorodnitsky1997}, le CBF peut se définir à partir de sa représentation harmonisable : \begin{equation}\label{eq:fbf}
B^H(\x)=\int_{\R^2}\frac{e^{i \x\cdot \xib}-1}{\|\xib\|^{H+1}}d\widehat W(\xib),\end{equation}
où $d\widehat W$ est une mesure brownienne complexe et $\x\cdot \xib$ désigne le produit scalaire sur $\R^2$. L'indice de Hurst $H$ est un paramètre fondamental du CBF, comme indicateur de rugosité de la texture. Plus $H$ est grand, plus la texture résultante est lisse. \\

Dans le but d'introduire de l'anisotropie dans le précédent modèle, Bonami et Estrade \cite{Bonami2003} ont remplacé l'indice de Hurst $H$ dans \eqref{eq:fbf} par une fonction dépendant de la direction $\xib$ et ont dérivé une nouvelle classe de \textit{champs browniens fractionnaires anisotropes} (CBFA) : 
\begin{equation}\label{eq:afbf}X(\x)=\int_{\R^2}\frac{e^{i \x\cdot \xib}-1}{\|\xib\|^{h(\arg \xib)+1}}d\widehat W(\xib).\end{equation}
Une classe de modèles anisotropes plus large peut être définie par
\begin{equation}\label{eq:afbfbis}
X(\x)=\int_{\R^2}(e^{i \x\cdot \xib}-1)f^{1/2}(\xib)~d\widehat W(\xib),
\end{equation}
où la densité spectrale $f$ est de la forme
\begin{equation}\label{eq:density}
f^{1/2}(\xib)=c(\arg \xib)\|\xib\|^{-h(\arg \xib)-1}.
\end{equation}
avec $c$ et $h$ deux fonctions $\pi$-périodiques définies sur l'intervalle $(-\pi/2,\pi/2]$ dont les images satisfont $c((-\pi/2,\pi/2])\subset \R^+$ et $h((-\pi/2,\pi/2])\subset (0,1)$. Quand $c$ et $h$ sont toutes deux constantes, on retrouve le CBF d'ordre $H\equiv h$.

\begin{figure}[t]
\begin{minipage}[b]{.48\linewidth}
  \centering
  \centerline{\includegraphics[width=4.0cm]{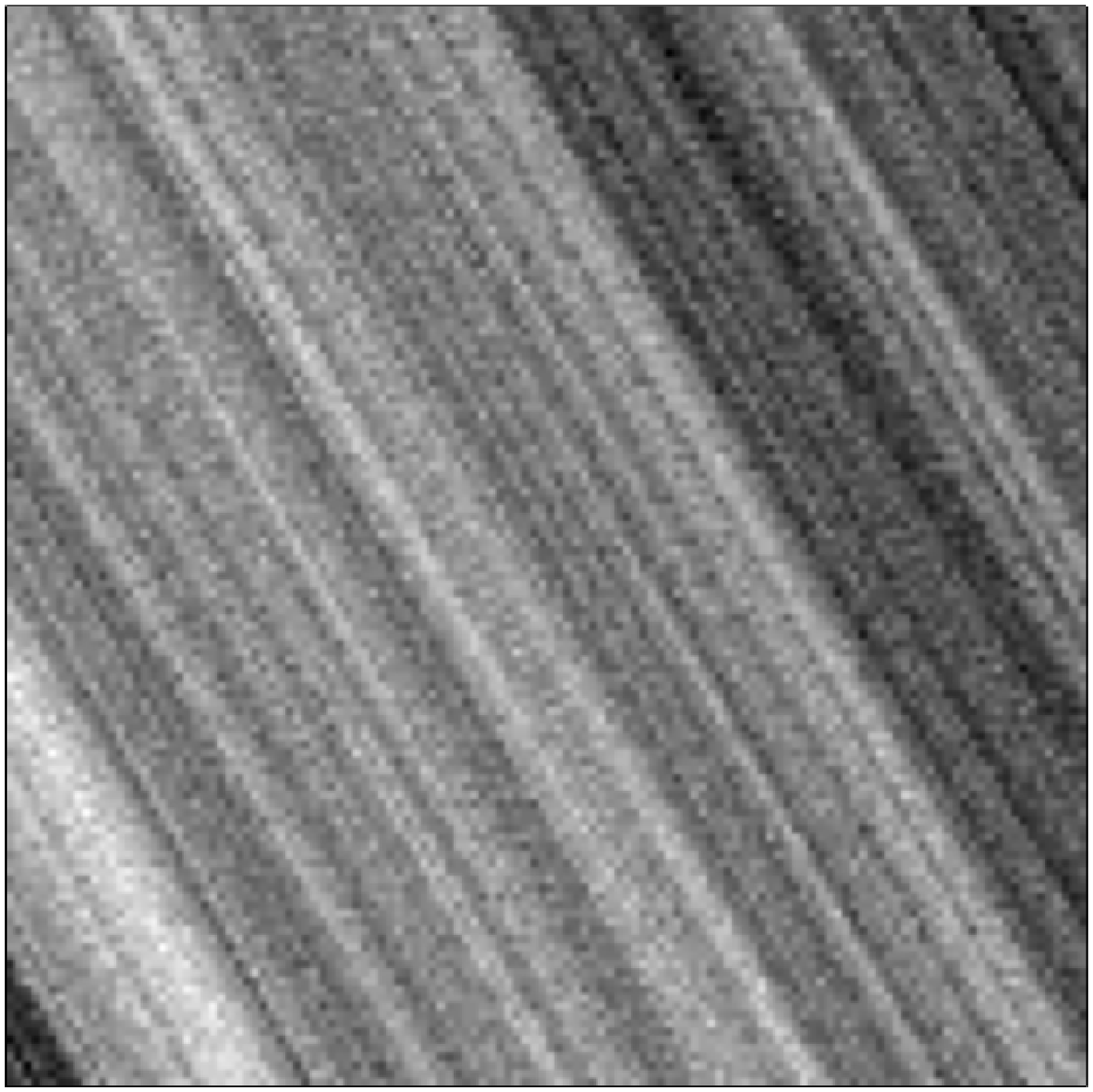}}
  \centerline{(a)}\medskip
\end{minipage}
\hfill
\begin{minipage}[b]{0.48\linewidth}
  \centering
  \centerline{\includegraphics[width=4.0cm]{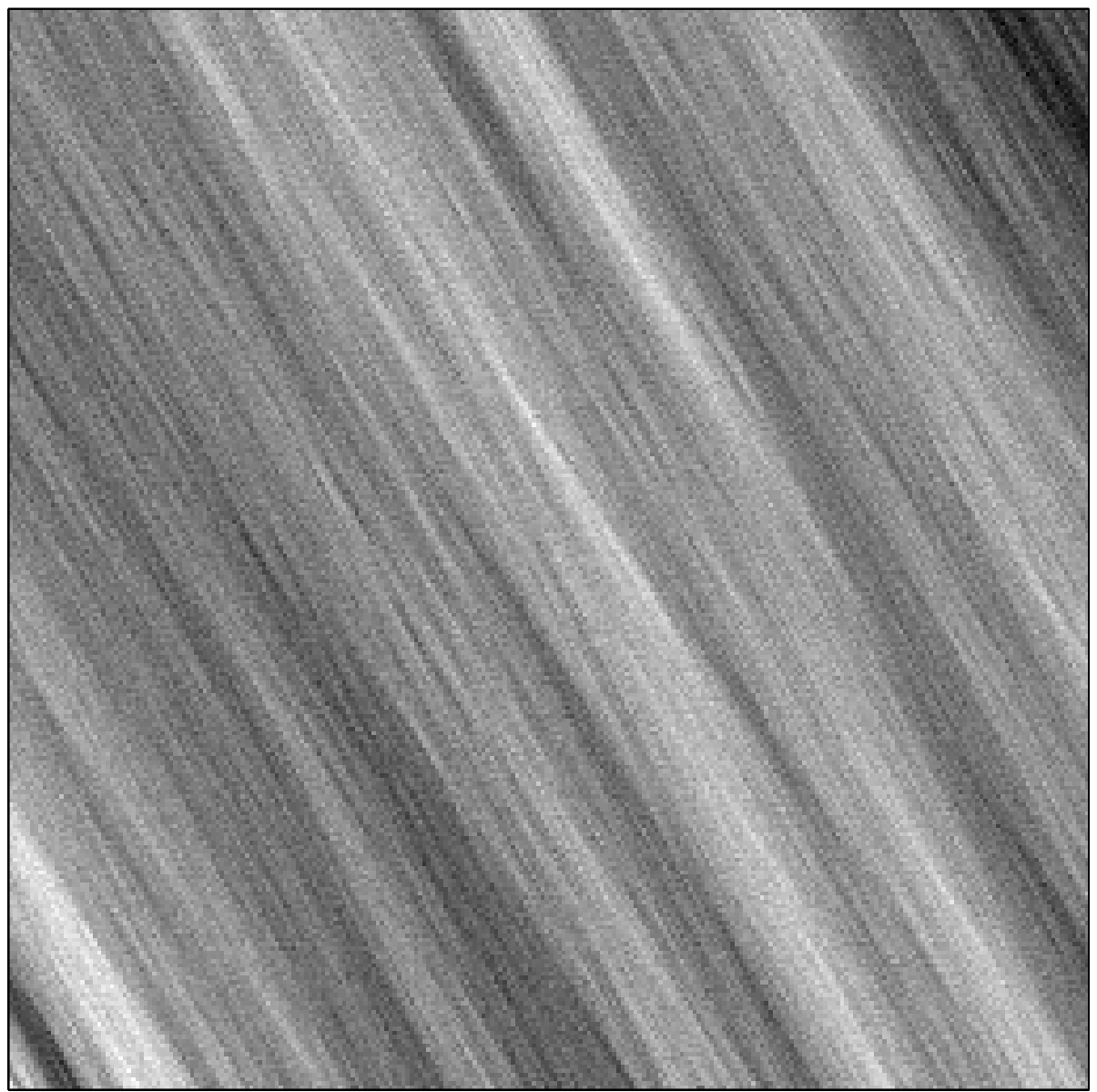}}
  \centerline{(b)}\medskip
\end{minipage}
\legende{Champ élémentaire de paramètres $H=0.2$, $\alpha_0=\frac \pi 6$ et $\alpha=10^{-2}$ simulé par (a) la méthode de Cholesky, (b) la méthode des bandes tournantes. }
\label{fig:fbf}
\end{figure}

Pour définir des modèles {\it  anisotropes stationnaires} possédant une orientation \textit{globale} $\alpha_0$, on peut fixer $h\equiv H$ dans~(\ref{eq:density}) et considérer :
 \begin{equation}\label{eq:elem}
 c_{\alpha_0,\alpha}(\arg(\xib))=\mathds{1}_{[-\alpha,\alpha]}(\arg(\xib)-\alpha_0),
 \end{equation}
 pour un certain $0<\alpha\leqslant \pi/2$.
 Nous retrouvons ainsi les \textit{champs élémentaires} de \cite{Bierme2012}, qui sont un cas particulier de CBFA. Quand $\alpha=\pi/2$, ce modèle correspond au CBF isotrope usuel d'indice de Hurst $H$, mais dès que $0<\alpha<\pi/2$, le champ n'est plus isotrope dans la mesure où les fréquences non nulles sont confinées dans le cône $[-\alpha+\alpha_0,\alpha+\alpha_0]$. 
 Plus la largeur du cône $\alpha$ décroît vers $0$, plus les fréquences sont concentrées suivant l'axe dirigé par $\alpha_0$, et plus la texture obtenue apparaît orientée dans la direction orthogonale à $\alpha_0$, conséquence des propriétés de la transformée de Fourier.\\
 Un exemple de champ d'orientation globale $\alpha_0=\frac \pi 6$ est donné figure~\ref{fig:fbf}
 (a) dont la réalisation est obtenue à partir de la fonction de covariance, qui est couramment employée pour simuler des champs gaussiens \cite{brouste2007}. Cependant la complexité élevée de ces algorithmes est un réel problème pour produire de grandes images de textures. Une méthode efficace a été proposée dans \cite{Bierme2012}, appelée méthode des bandes tournantes, et est utilisée ici pour simuler la texture de la figure~\ref{fig:fbf} (b).
 
%

\begin{figure*}[t]
  \centering
  \includegraphics[scale=0.7]{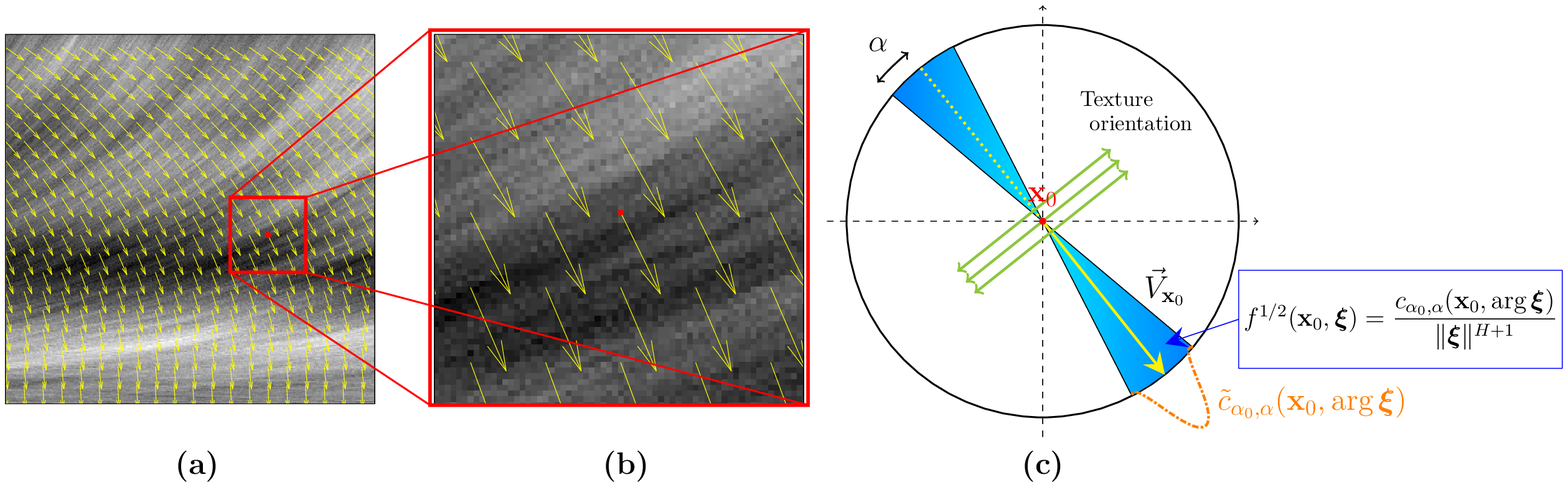}
  \legende{(a) Texture obtenue à partir du champ d'orientations $\vec{V}_{(x,y)}^1$ en jaune et $H=0.2$, (b) zoom autour du point rouge $\x_0=(x,y)$ montrant localement un champ élémentaire orienté, (c) diagramme illustrant chaque paramètre du modèle CBFLA.}
  \label{fig:summary}
\end{figure*}

\section{Une nouvelle classe de champs gaussiens à orientation prescrite}
\label{sec:LAFBF}


On introduit maintenant un nouveau modèle gaussien comme une adaptation locale du champ élémentaire défini dans~(\ref{eq:afbfbis})-(\ref{eq:density})-(\ref{eq:elem}) avec $h\equiv H$. Plus précisément, on définit le {\it champ brownien fractionnaire localement anisotrope} (CBFLA) par :
\begin{equation}\label{eq:LAFBF}
X(\x)=\int_{\R^2}(e^{i\x\cdot\xib}-1)f^{1/2}(\x,\xib)~d\widehat{W}(\xib),
\end{equation}
avec \begin{equation}\label{eq:dens}f^{1/2}(\x,\xib)=c_{\alpha_0,\alpha}(\x,\arg \xib)\|\xib\|^{-H-1},\end{equation}
\begin{equation}\label{eq:cfunc}c_{\alpha_0,\alpha}(\x,\arg \xib)=\mathds{1}_{[-\alpha,\alpha]}(\arg(\xib)-\alpha_0(\x)),\end{equation} $\alpha_0$ étant maintenant une fonction différentiable sur $\R^2$.
Le CBFLA est dérivé du CBFA de la même façon que le \textit{mouvement brownien multifractionnaire} (MBM) est dérivé du \textit{mouvement brownien fractionnaire} (MBF), puisque le paramètre d'orientation $\alpha_0$ est remplacé par une fonction dépendant de la position $\x$, tout comme l'indice de Hurst $H$ du MBM peut varier suivant la position \cite{peltier1995}.\\

\textbf{Champ tangent en chaque point}. Afin de décrire et simuler cette nouvelle classe de champs gaussiens, nous avons besoin de la notion de champ tangent, que l'on rappelle brièvement : un champ aléatoire $X$ est dit localement asymptotiquement autosimilaire d'ordre $H\in (0,1)$ en $\x_0\in \R^2$ si pour tout $\h\in\R^2$ le champ aléatoire $\rho^{-H}(X(\x_0+\rho \h)-X(\x_0))$ admet en loi une limite $Y_{\x_0}$ quand $\rho\to 0$ (c.f \cite{benassi1997,falconer2002} pour une définition plus générale). Le champ $Y_{\x_0}$ est alors appelé le champ tangent de $X$ en $\x_0$. Ainsi, le champ aléatoire $X$ se \textit{comporte localement comme} $Y_{\x_0}$ quand $\x\to \x_0$. Cette notion a été introduite pour la première fois dans \cite{benassi1997} pour décrire le comportement local du mouvement brownien multifractionnaire (qui se comporte localement comme un MBF).

Nous avons prouvé dans \cite{polisano:2015} que le CBFLA $X$ de \eqref{eq:LAFBF} admet pour champ tangent :
\begin{equation}\label{eq:tangent}
Y_{\x_0}(\x)=\displaystyle \int_{\R^2}(e^{i \x\xib}-1)f^{1/2}(\x_0,\xib)~d\widehat{W}(\xib)\;.
\end{equation}
On remarque que ce champ tangent $Y_{\x_0}$ est un champ élémentaire (\ref{eq:afbfbis})-(\ref{eq:density})-(\ref{eq:elem}). La simulation d'un CBFLA requiert maintenant de savoir simuler le champ tangent en tout point d'une grille discrète $r^{-1}\mathbb{Z}^2\cap [0,1]^2$ avec $r=2^k-1, k\in \mathbb{N}^{\star}$. \\
Deux méthodes sont proposées : \\ 


\textbf{ (i) Simulation des champs tangents par la méthode de Cholesky}. La covariance d'un champ élémentaire $Y_{\x_0}$, qui est à accroissements stationnaires, est reliée à son variogramme  $v_{Y_{\x_0}}$ par $\text{Cov}(Y_{\x_0}(\x),Y_{\x_0}(\y))=v_{Y_{\x_0}}(\x)+v_{Y_{\x_0}}(\y)-v_{Y_{\x_0}}(\x-\y)$, avec $v_{Y_{\x_0}}$ déterminé par une formule explicite \cite{Bierme2012}. La matrice de covariance $\Sigma$ associée est de taille $r^2\times r^2$, et sa décomposition de Cholesky $\Sigma=LL^T$ est de complexité $O(r^6)$, ce qui est extrêmement coûteux. Une réalisation du champ élémentaire peut être obtenue par $Y_{\x_0}\sim L\mathbf{Z}$ (figure \ref{fig:fbf} (a)), où $\mathbf{Z}\sim (\mathcal{N}(0,1))^{r^2}$. Une méthode plus fine est proposée dans  \cite{brouste2007}. \\

\textbf{(ii) Simulation des champs tangents par bandes tournantes}. La méthodologie de \cite{Bierme2012} a été adaptée dans \cite{polisano:2014}.

\textit{-- Formulation discrète du champ tangent}\\
On peut dériver une expression intégrale du variogramme de $Y_{\x_0}$ par un changement de variable en polaire : 
\begin{equation}\label{eq:vario}
\begin{array}{rcl}v_{Y_{\x_0}}(\x)&=&\displaystyle \frac{1}{2}\int_{\R^2}|e^{i \x \cdot \xib}-1|^2f(\x_0,\xib)d\xib \\ &=& \displaystyle \gamma(H)\int_{-\pi/2}^{\pi/2}c_{\alpha_0,\alpha}(\x_0,\theta)~|\x \cdot \ub(\theta)|^{2H} d\theta\end{array},
\end{equation}
où $\ub(\theta)=(\cos \theta,\sin \theta)$ et $\gamma(H)=\frac{\pi}{2H\Gamma(2H)\sin(H\pi)}$. \\
L'intégrale (\ref{eq:vario}) est de la forme $\int_{-\pi/2}^{\pi/2}\tilde v_{\theta}(\x\cdot \ub(\theta))d\theta$ avec\\
$\tilde v_{\theta}=\gamma(H)c_{\alpha_0,\alpha}(\x_0,\theta)|\cdot |^{2H}$. En ignorant le facteur \\
$\gamma(H)c_{\alpha_0,\alpha}(\x_0,\theta)$, on reconnait que $\tilde v_{\theta}$ est le variogramme d'un CBF d'ordre $H$.  Par conséquent, $Y_{\x_0}$ peut être vu comme une somme de CBF tournant autour de l'origine. En discrétisant $\theta$ en un ensemble ordonné $(\theta_i)_{1\leqslant i\leqslant n}$ de $n$ bandes d'orientation, et en notant $(\lambda_i)_{1\leqslant i\leqslant n}$ les largeurs de bandes associées $\lambda_i=\theta_{i+1}-\theta_i$, les {\it champs à bandes tournantes} s'écrivent
\vspace{-0.2cm}
 \begin{equation}\label{eq:turning}Y_{\x_0}(\x)=\gamma(H)^{\frac{1}{2}}\sum_{i=1}^n\sqrt{\lambda_ic_{\alpha_0,\alpha}(\x_0,\theta_i)}B_i^{H}(\x\cdot \ub(\theta_i)),
 \end{equation} 
où les $B_i^{H}$'s sont $n$ CBF d'ordre $H$. Cette version discrète est une bonne approximation dès lors que $\displaystyle \max_{i}\lambda_i \leqslant \varepsilon$ petit.

\begin{algorithm}[t]
 \caption{Simulation du CBFLA}
 \label{alg:cbfla}
 \begin{algorithmic}[1]
 \INPUT $r=2^k-1$, $H$, $\alpha_0$, $\alpha$, $\epsilon$
 \OUTPUT $X$ CBFLA de taille $(r+1)\times (r+1)$
   \STATE $(p_i,q_i)_{1\leqslant i\leqslant n} \leftarrow \text{ChoixBandesDynamique}(r,\varepsilon)$
   \STATE Calculer et trier les angles $(\theta_i)_{1\leqslant i\leqslant n}$ : $\theta_i\leftarrow \text{atan2}(p_i,q_i)$
   \STATE Calculer les largeurs de bandes $(\lambda_i)_{1\leqslant i\leqslant n}$ : $\lambda_i\leftarrow \theta_{i+1}-\theta_i$
   \STATE Générer $n$ CBF : $B_i^H\leftarrow \text{circMBF}(r(|p_i|+|q_i|),H)$
   \STATE Initialisation : $X\leftarrow 0$
   \FORALL {$(k_1,k_2)$}
		\FOR{$i=1$ \TO $n$}
		        \STATE $\omega_i\leftarrow \sqrt{\lambda_i\gamma(H)c_{\alpha_0,\alpha}((k_1,k_2),\theta_i)}\left(\frac{\cos \theta_i}{rq_i}\right)^H$
   			\STATE $X(k_1,k_2)\leftarrow X(k_1,k_2)+\omega_iB_i^H(k_1q_i+k_2p_i)$
		 \ENDFOR
   \ENDFOR
 \end{algorithmic}
 \end{algorithm}

\textit{-- Simulation le long de bandes particulières}\\
Suivant \cite{Bierme2012}, on choisit $\theta_i$ tel que $\tan(\theta_i)=\frac{p_i}{q_i}$, avec $p_i,q_i\in \mathbb{Z}$,  \begin{equation}\label{eq:trick}\begin{array}{c}\displaystyle \left\{B_i^{H}\left(\frac{k_1}{r}\cos \theta_i+\frac{k_2}{r}\sin \theta_i\right);0\leqslant k_1,k_2\leqslant r\right\}\overset{\mathcal{L}}{=}\\ \displaystyle
 \left(\frac{\cos \theta_i}{rq_i}\right)^{H}\{B_i^{H}(k_1q_i+k_2p_i);0\leqslant k_1,k_2\leqslant r\}\end{array},\end{equation} et $B_i^H(\x\cdot \ub(\theta_i))$ peut ainsi être généré en utilisant l'algorithme rapide de Perrin \textit{et al.} \cite{perrin2002} sur une grille régulière.

\textit{-- Choix dynamique des bandes discrètes}\\
En pratique, le choix des bandes d'orientation $(\theta_i)_{1\leqslant i\leqslant n}$ est gouverné par le coût de calcul global des $B_i^H$, au travers de la programmation dynamique  \cite{Bierme2012}.\medskip

 \begin{figure}[t]

  \centering
  \includegraphics[scale=0.34]{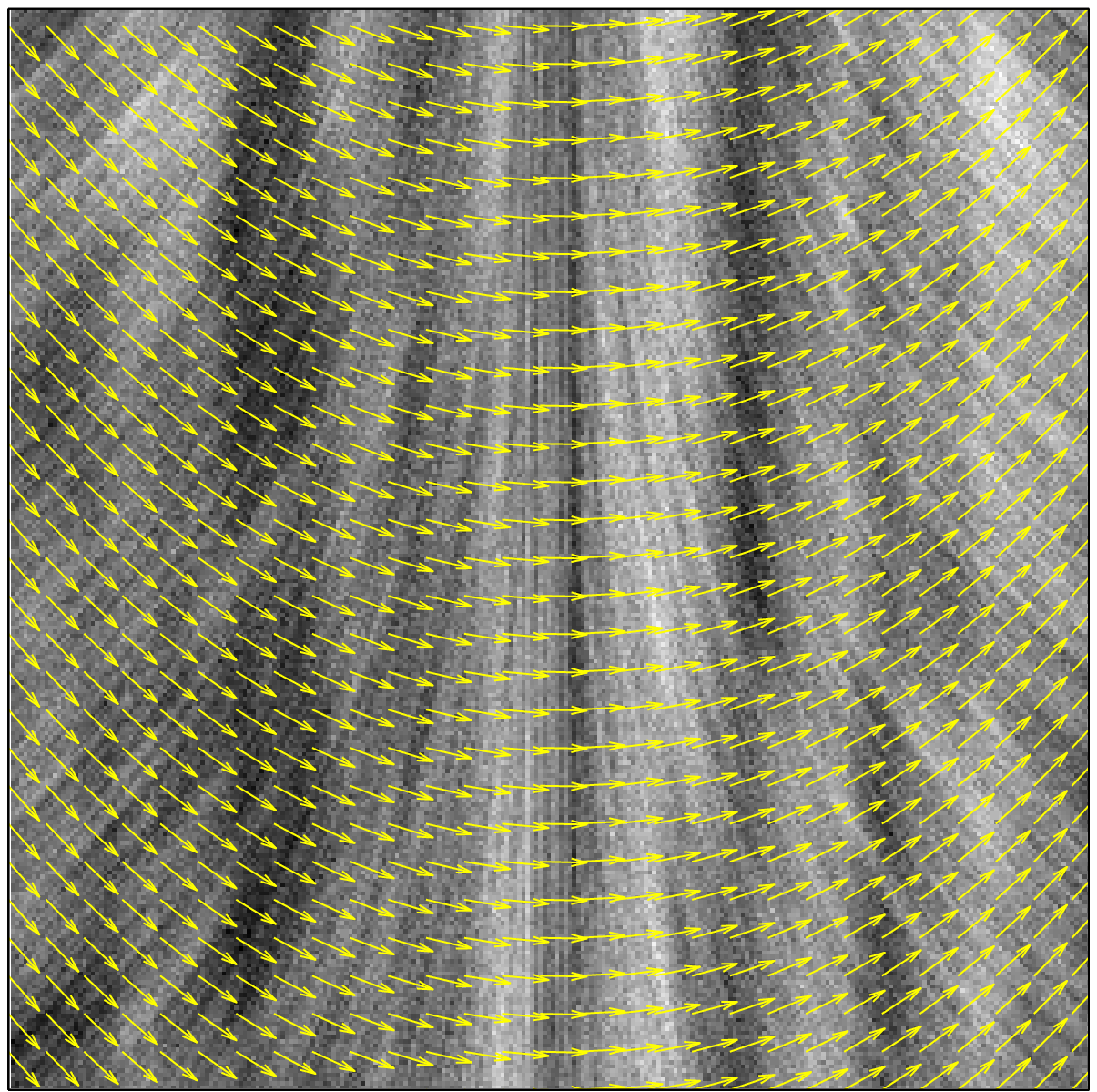}
\caption{Texture obtenue pour le champ $\vec{\V}_{(x,y)}^2$ et $H=0.2$.}
\label{fig:vague}
\end{figure}

\textbf{Simulation du CBFLA}. Comme observé dans \cite{peltier1995} pour le MBM, un champ gaussien peut être simulé à partir de ses champs tangents. Le CBFLA se comportant localement comme son champ tangent, en chaque pixel $\x_0$, on affecte $X(\x_0)=Y_{\x_0}(\x_0)$ obtenu par la méthode de Cholesky (i) ou des bandes tournantes (ii). Le pseudo-code est donné dans l'algorithme \ref{alg:cbfla}. L'étape de pré-traitement (instructions 1,2,3,4 ), laquelle ne dépend pas des orientations locales prescrites, est exécutée une fois pour toutes. Le reste de l'algorithme a une complexité linéaire en $O(r^2\log n)$. En effet, en chaque point $(k_1,k_2)$, une bande tournante $\theta_i$ intervient dans le calcul de $X(k_1,k_2)$ si et seulement si elle appartient au cône $c_{\alpha_0,\alpha}((k_1,k_2),\theta_i)\neq 0$, soit $|\theta_i-\alpha_0((k_1,k_2))|\leqslant \alpha$. Le tableau $\theta_i$ étant trié, un tel indice $i$ peut être identifié par une recherche dichotomique. 


\begin{figure}[t]
\begin{minipage}[b]{.48\linewidth}
  \centering
  \centerline{\includegraphics[width=4.0cm]{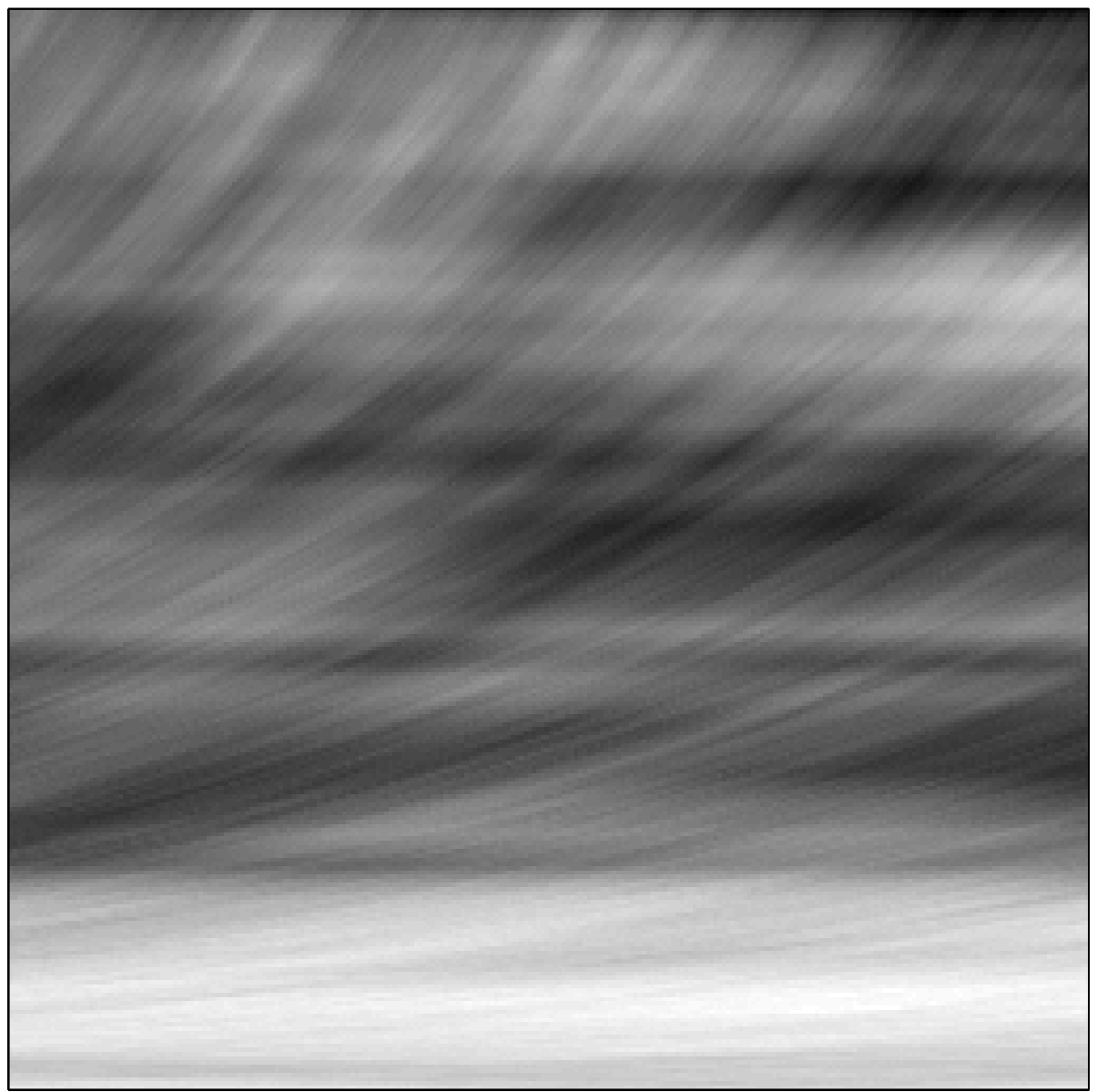}}
  \centerline{(a)}\medskip
\end{minipage}
\hfill
\begin{minipage}[b]{0.48\linewidth}
  \centering
  \centerline{\includegraphics[trim = 0mm 0mm 30mm 30mm, clip, width=4.0cm]{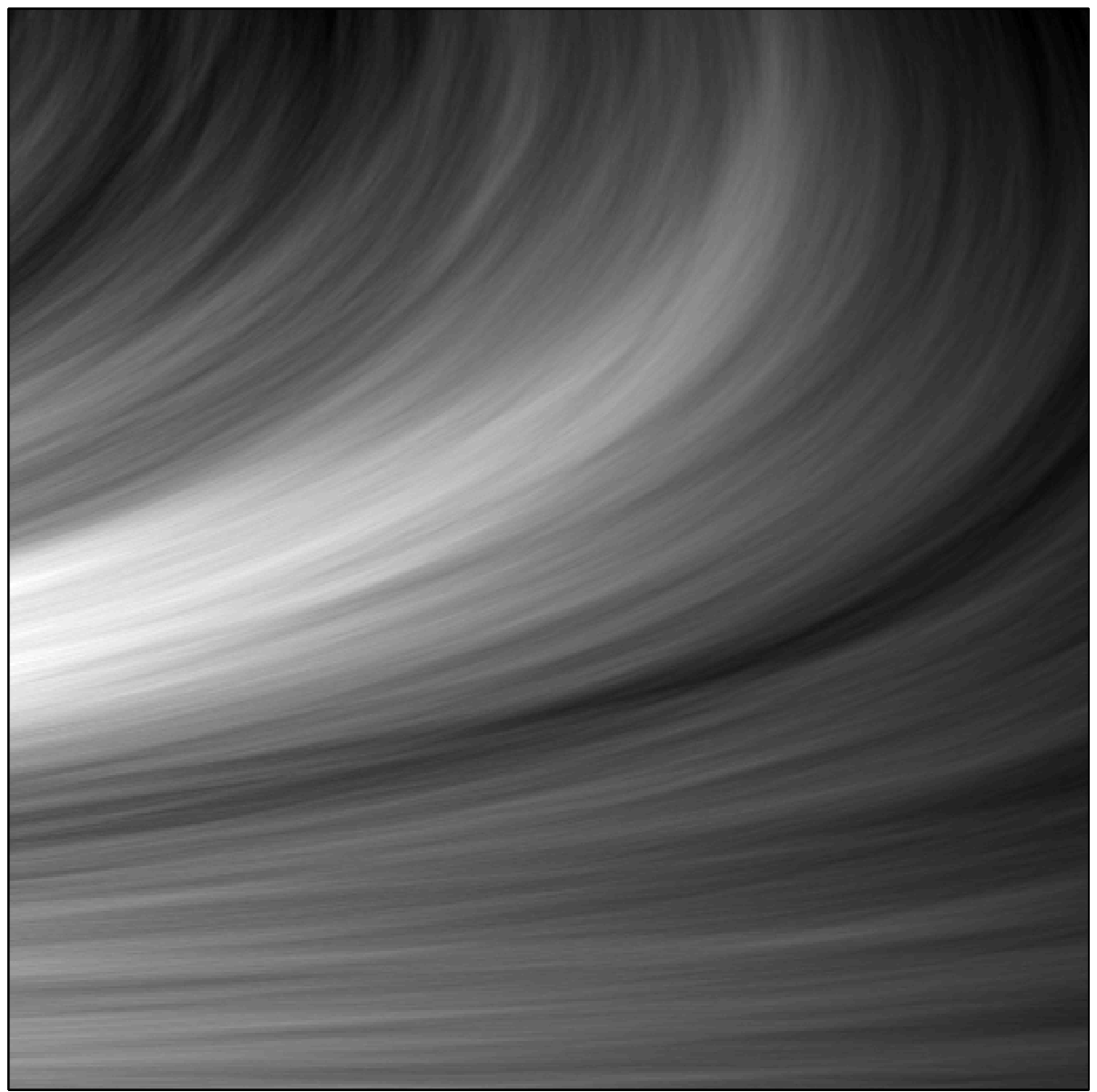}}
  \centerline{(b)}\medskip
\end{minipage}
\legende{Texture obtenue pour le champ $\vec{\V}_{(x,y)}^1$ et $H=0.7$ (a) par bandes tournantes : bandes d'artefacts (b) par Cholesky}
\label{fig:turningVSchol}
\end{figure}


\section{Synthèse de textures orientées}

Les paramètres utilisés dans les simulations sont $r=255$, $\alpha=10^{-1}$, et $\varepsilon=10^{-2}$. Afin d'éviter les artefacts numériques dans (\ref{eq:turning}) on considère une version régularisée $\tilde c_{\alpha_0,\alpha}$ de la fonction indicatrice $c_{\alpha_0,\alpha}$, typiquement une gaussienne. Pour $\alpha_0$ constant, on retrouve les résultats de \cite{Bierme2012} (figure~\ref{fig:fbf} (b)), qui est conforme à ce que l'on obtient par la méthode de Cholesky (figure~\ref{fig:fbf} (a)) considérée comme exacte.  Toutefois l'emploi d'une fenêtre gaussienne $\tilde c_{\alpha_0,\alpha}$ lisse la texture et augmente la régularité. Nous présentons maintenant des réalisations de textures à orientation prescrite en chaque point $\x_0$, donnée par le champ de vecteurs $\vec{\V}_{\x_0}=\ub(\alpha_0(\x_0))$. La figure~\ref{fig:summary} (a) est la texture obtenue à partir du champ de vecteurs $\vec{\V}_{(x,y)}^1=(\cos(-\pi/2+y),\sin(-\pi/2+y))$. Un zoom autour du point  $\x_0$ (figure~\ref{fig:summary} (b)) montre que localement un CBFLA se comporte comme un champ élémentaire. La figure~\ref{fig:summary} (c) schématise la fonction de densité locale au $\x_0$ et les différents paramètres.  Le champ $\vec{\V}_{(x,y)}^1$ présente des bandes d'artefacts plus marquées pour des indices de Hurst supérieurs ($H=0.7$), comme le montre la figure~\ref{fig:turningVSchol} (a) tandis que la simulation par la méthode de Cholesky en est exempte figure~\ref{fig:turningVSchol} (b). On considère enfin le champ de vecteurs, $\vec{\V}_{(x,y)}^2=(\cos(\sin(2x-1)),\sin(\sin(2x-1))$, dont la texture correspondante est représentée figure~\ref{fig:vague}. D'autres exemples ont été donnés dans \cite{polisano:2014}. Comme attendu, les textures obtenues avec notre approche présentent un comportement local anisotrope, orienté orthogonalement au champ de vecteurs. De plus, la simulation d'une texture de taille $256\times 256$ prend seulement quelques secondes avec la méthode des bandes tournantes, et plusieurs heures avec la méthode de Cholesky.

\section{Conclusion}

Nous avons introduit un nouveau modèle stochastique pour la simulation de textures à orientation locale prescrite, vu comme une version locale du champ brownien fractionnaire anisotrope. Nous avons tiré parti de la formulation par champ tangent couplée à la méthode des bandes tournantes, pour proposer un algorithme efficace de simulation de ces textures. En couplant avec la méthode de Cholesky on supprime les artefacts mais on augmente fortement le temps de calcul. Une extension de notre modèle à des indices de Hurst variant localement est en cours d'étude \cite{polisano:2015}, et servira de test pour des estimateurs de régularité locale \cite{meyer1996wavelet} et d'orientation locale \cite{unser2009multiresolution,olhede2014detecting}.\\

\noindent \textit{Remerciements} -- Ce travail est soutenu par l'ANR dans le cadre du projet ANR-13-BS03-0002-01 (ASTRES).

\NoAutoSpaceBeforeFDP
\renewcommand{\baselinestretch}{0}
\scalefont{0.95}
\bibliographystyle{IEEEbib}
\bibliography{refs}

%
%
%

\end{document}